\documentclass[10pt]{article}
\usepackage{latexsym}
\usepackage{amsfonts}
\usepackage{enumerate}
\usepackage{multicol}
\usepackage{graphicx}
\usepackage{amssymb}
\usepackage{amsmath}
\usepackage{epic}

\title{On the minimum size of restricted sumsets in cyclic groups \\[.4in]}

\author{B\'{e}la Bajnok \\[.1in] Department of Mathematics, Gettysburg College \\
Gettysburg, PA 17325-1486 USA \\E-mail:  bbajnok@gettysburg.edu \\[.4in]}

\date{May 13, 2013}

\newtheorem{thm}{Theorem}

\newtheorem{lem}[thm]{Lemma}
\newtheorem{cor}[thm]{Corollary}
\newtheorem{prop}[thm]{Proposition}
\newtheorem{conj}[thm]{Conjecture}

\newtheorem{prob}[thm]{Problem}

\begin{document}

\maketitle

\begin{abstract}

For positive integers $n$, $m$, and $h$, we let $\rho \hat{\;}(\mathbb{Z}_n, m, h)$ denote the minimum size of the $h$-fold restricted sumset among all $m$-subsets of the cyclic group of order $n$.  The value of $\rho \hat{\;}(\mathbb{Z}_n, m, h)$ was conjectured for prime values of $n$ and $h=2$ by Erd\H{o}s and Heilbronn in the 1960s; Dias da Silva and Hamidoune proved the conjecture in 1994 and generalized it for an arbitrary $h$, but little is known about the case when $n$ is composite.  Here we exhibit an explicit upper bound for all $n$, $m$, and $h$; our bound is tight for all known cases (including all $n$, $m$, and $h$ with $n \leq 40$).  We also provide counterexamples for conjectures made by Plagne and by Hamidoune, Llad\'o, and Serra.

\end{abstract}



\thispagestyle{empty}
\pagestyle{myheadings}
\markright{Minimum size of restricted sumsets}

\section{Introduction}

Let $G$ be a finite abelian group written with additive notation.  For a positive integer $h$ and a nonempty subset $A$ of $G$, we let $hA$ and $h\hat{\;}A$ denote the $h$-fold unrestricted sumset and the $h$-fold restricted sumset of $A$, respectively; that is, $hA$ is the collection of sums of $h$ not-necessarily-distinct elements of $A$, and $h\hat{\;}A$ consists of all sums of $h$ distinct elements of $A$.  For a positive integer $m \leq |G|$, we 
let $$\rho(G, m, h) = \min \{ |hA| \; : \; A \subseteq G, |A|=m\}$$ and   
$$\rho \hat{\;}(G, m, h) = \min \{ |h \hat{\;}A|  \; : \; A \subseteq G, |A|=m\}$$ (as usual, $|S|$ denotes the size of the finite set $S$).

The value of $\rho(G, m, h)$ has a long and distinguished history and has been determined for all $G$, $m$, and $h$, but the quantity $\rho \hat{\;}(G, m, h)$ remains largely unknown to this day.  In this paper we attempt to find good upper bounds for $\rho \hat{\;}(G, m, h)$ in the case when $G$ is cyclic.   

We start by a brief review of the case of unrestricted sumsets.  In 1813, Cauchy  \cite{Cau:1813a} found the value of $ \rho (\mathbb{Z}_p, m, 2)$ for prime $p$; in 1935 Davenport \cite{Dav:1935a} rediscovered Cauchy's result, which is now known as the Cauchy--Davenport Theorem.  (Davenport was unaware of Cauchy's work until twelve years later; see \cite{Dav:1947a}.)  

\begin{thm}[Cauchy--Davenport Theorem] \label{Cauchy--Davenport}
For a prime $p$ and any positive integer $m \leq p$ we have $$ \rho (\mathbb{Z}_p, m, 2)=\min\{p,2m-1\}.$$
\end{thm}

After various partial results, the general case was finally solved in 2003 by Eliahou, Kervaire, and Plagne \cite{EliKerPla:2003a} (see also \cite{Pla:2003a}, \cite{EliKer:2007a}, and  \cite{Pla:2006a}).  To state the result, we introduce the function
$$u(n,m,h)=\min \{f_d \; : \; d \in D(n)\},$$ where $n$, $m$, and $h$ are positive integers, $D(n)$ is the set of positive divisors of $n$, and
$$f_d=\left(h\left \lceil \tfrac{m}{d} \right \rceil-h +1 \right) \cdot d.$$  (Here $u(n,m,h)$ is a relative of the Hopf--Stiefel function used also in topology and bilinear algebra; see, for example, \cite{Sha:1984a}, \cite{Pla:2003a}, and \cite{Kar:2006a}.)

\begin{thm} [Eliahou, Kervaire, and Plagne; cf.~\cite{EliKerPla:2003a}] \label{value of u}
Let $n$, $m$, and $h$ be positive integers with $m \leq n$.  For any abelian group $G$ of order $n$ we have
$$\rho (G, m, h)=u(n,m,h).$$
\end{thm}  

We can observe that we have $f_1=hm-h+1$ and $f_n=n$; therefore, Theorem \ref{value of u} is indeed a generalization of the Cauchy--Davenport Theorem.

The fact that 
\begin{eqnarray} \label{rho>=u}
\rho (G, m, h) \geq u(n,m,h)
\end{eqnarray} follows readily from the following version of Kneser's Theorem (see \cite{Kne:1956a} and \cite{Pla:2006a}).  (As customary, for subsets $X$ and $Y$ of $G$ we let $X+Y$ denote the set of sums $x+y$ with $x \in X$ and $y \in Y$; if $X=\{x\}$, we also write $x+Y$ instead of $X+Y$.  The stabilizer of a subset $X$ of $G$ is the largest subset---as it turns out, subgroup---$H$ of $G$ for which $X+H=X$.)  

\begin{thm}  [Generalized Kneser's Theorem; cf.~\cite{Pla:2006a}] \label{Kneser}
Let $G$ be an abelian group, $m$ and $h$ positive integers, and $A$ be an $m$-subset of $G$ for which
$|hA| \leq hm-h.$
Then the stabilizer $H$ of $hA$ is a nontrivial subgroup of $G$, and, with $A/H$ denoting the image of $A$ under the canonical projection $G \rightarrow G/H$, we have
$$|h(A/H)| \geq h |A/H|-h+1.$$
\end{thm}

To see that Theorem \ref{Kneser} implies (\ref{rho>=u}), choose an $m$-subset $A$ of $G$ for which $\rho (G, m, h)=|hA|$.  Clearly, if $|hA| \geq  hm-h+1=f_1$, then    
$\rho (G, m, h)=|hA|\geq u(n,m,h)$, so assume that $|hA| \leq hm-h$.  In this case, we can apply Theorem \ref{Kneser}: letting $d$ denote the order of $H$ and noting that $A$ intersects at least $\left \lceil m/d \right \rceil$ cosets of $H$ and thus $|A/H| \geq \left \lceil m/d \right \rceil$, we have 
$$|hA| = |hA+H|=|h(A/H)| \cdot d \geq (h |A/H|-h+1) \cdot d \geq f_d \geq u(n,m,h).$$

Therefore, to establish Theorem \ref{value of u}, we need to prove the existence of an $m$-subset $A$ of $G$ with $|hA|=u(n,m,h)$; to do so, we may proceed as follows.  For a given divisor $d$ of $n$, we write $m$ as $m=cd+k$ with $1 \leq k \leq d$, and set
$$A_d(n,m)=\bigcup_{i=0}^{c-1} (i+H) \cup \left\{c + j \cdot \tfrac{n}{d}  \; : \;  j=0,1,2,\dots,k-1 \right\},$$ where $H$ is the (unique) subgroup of order $d$ of $\mathbb{Z}_n$.  As we verify in Section  \ref{regular sets}, we have  
\begin{eqnarray} \label{|hA_d(n,m)|}
 |hA_d(n,m)| &=& \min\{n,f_d,hm-h+1\};
\end{eqnarray}
since $f_n=n$ and $f_1=hm-h+1$, we get  
$$\rho(\mathbb{Z}_n, m, h) \leq \min \{f_d  \; : \;  d \in D(n)\}=u(n,m,h).$$ 
This, with (\ref{rho>=u}), establishes Theorem \ref{value of u} for cyclic groups.

In this paper we find $|h\hat{\;}A_d(n,m)|$, the size of the $h$-fold restricted sumset of $A_d(n,m)$; as it turns out, it will yield a very good upper bound (though sometimes not the exact value) for $\rho \hat{\;}(\mathbb{Z}_n, m, h)$.  Our computations result in the following.

\begin{thm} \label{sizeofh^A_d}
Suppose that $n$, $m$, and $h$ are positive integers, $m \leq n$, and $d \in D(n)$; let $A_d(n,m)$ be defined as above.  

If $h=m$, then $|h\hat{\;}A_d(n,m)| =1$, and if $h>m$, then $|h\hat{\;}A_d(n,m)| =0$.  

For  $h \leq m-1$, let $k$ and $r$ denote the positive remainder of $m$ mod $d$ and $h$ mod $d$, respectively; we then have  
$$|h\hat{\;}A_d(n,m)| = \left\{\begin{array}{ll}
\min\{n, f_d, hm-h^2+1\}  & \mbox{if $h \leq \min\{k,d-1\}$;}\\ \\
\min\{n, hm-h^2+1 - \delta_d\} & \mbox{otherwise;}\\  
\end{array}\right.$$
where $\delta_d$ is the ``correction term'' defined as 
$$\delta_d= \left\{ 
\begin{array}{ll}
(k-r)r-(d-1)  & \mbox{if $r<k$,}\\ 
(d-r)(r-k)-(d-1)  & \mbox{if $k<r<d$,}\\
d-1  & \mbox{if $k=r=d$,}\\
0  & \mbox{otherwise.}  
\end{array}\right.$$   
\end{thm}
We prove Theorem \ref{sizeofh^A_d} in Section  \ref{regular sets}.

As an analogue of the function $u(n,m,h)$ defined above, we introduce the function
$$u \hat{\;}(n,m,h)=\min \{f \hat{_d} \; : \; d \in D(n)\},$$ where $n$, $m$, and $h$ are positive integers, $m \leq n$, $D(n)$ is the set of positive divisors of $n$, and
$$f \hat{_d}=|h\hat{\;}A_d(n,m)|$$ is given explicitly by Theorem \ref{sizeofh^A_d}.  

It is easy to see that $u \hat{\;}(n,m,1)=m$ for all $n$ and $m$, and we will find simple expressions for $u \hat{\;}(n,m,2)$ and $u \hat{\;}(n,m,3)$ in Section \ref{h=2 and h=3}.  We can also observe that $$f \hat{_1}=f \hat{_n}=\min\{n,hm-h^2+1\},$$ and thus for prime $p$ values we have 
\begin{eqnarray} \label{u-hat(p,m,h)}
u \hat{\;}(p,m,h) &=& \min\{p,hm-h^2+1\}.
\end{eqnarray} 

Obviously, $u \hat{\;}(n,m,h)$ provides an upper bound for $\rho \hat{\;}(\mathbb{Z}_n, m, h)$, and it seems that it is a remarkably good one.  In an exhaustive computer search for all $n \leq 40$, all $m$-subsets of $\mathbb{Z}_n$ (with $m \leq n$), and all $h$ (with $1 \leq h \leq m$), we found that $u \hat{\;}(n,m,h)$ agrees with $\rho \hat{\;}(\mathbb{Z}_n, m, h)$ in the overwhelming majority (over 99\%) of cases, and when it does not, we have $\rho \hat{\;}(\mathbb{Z}_n, m, h)=u \hat{\;}(n,m,h)-1$.  In light of this we pose the following: 

\begin{prob} \label{problem}
Classify all situations when $$\rho \hat{\;}(\mathbb{Z}_n, m, h)< u \hat{\;}(n,m,h);$$ in particular, how much smaller than $u \hat{\;}(n,m,h)$ can $\rho \hat{\;}(\mathbb{Z}_n, m, h)$ be?  

\end{prob}
As we just mentioned, in all cases that we are aware of, we have $\rho \hat{\;}(\mathbb{Z}_n, m, h)= u \hat{\;}(n,m,h)$ or $\rho \hat{\;}(\mathbb{Z}_n, m, h)= u \hat{\;}(n,m,h)-1$.
 
Let us mention three---as it turns out, quite predicative---examples for the case when $\rho \hat{\;}(\mathbb{Z}_n, m, h)= u \hat{\;}(n,m,h)-1$.
\begin{itemize}
\item $u \hat{\;}(12,7,2)=11$, but $\rho \hat{\;}(\mathbb{Z}_{12}, 7, 2)=10$ as shown by the set $$C_1=\{0,4\}\cup\{1,5,9\}\cup\{6,10\}.$$ 
\item $u \hat{\;}(10,6,3)=10$, but $\rho \hat{\;}(\mathbb{Z}_{10}, 6, 3)=9$ as shown by the set $$C_2=\{0,2,4,6\}\cup\{7,9\}.$$ 
\item $u \hat{\;}(15,8,3)=15$, but $\rho \hat{\;}(\mathbb{Z}_{15}, 8, 3)=14$ as shown by the set $$C_3=\{0,3,6,9\}\cup\{10,13,1,4\}.$$ 
\end{itemize}
Needless to say, we chose to represent our sets in the particular formats above for a reason.  In fact, all known situations where $\rho \hat{\;}(\mathbb{Z}_n, m, h)< u \hat{\;}(n,m,h)$ can be understood by a particular modification of our sets $A_d(n,m)$, as we now describe.

Observe that the $m$ elements of $A_d(n,m)$ are within $\lceil m/d \rceil$ cosets of the order $d$ subgroup $H$ of $\mathbb{Z}_n$, and at most one of these cosets does not lie entirely in $A_d(n,m)$.  We now consider the situation when the $m$ elements are still within $\lceil m/d \rceil$ cosets of $H$, but exactly two of these cosets don't lie entirely in our set.  In order to do so, we write $m$ in the form
$$m=k_1+(c-1)d+k_2$$ for some positive integers $c$, $k_1$, and $k_2$; we assume that $k_1<d$, $k_2<d$, but $k_1+k_2>d$.  We then are considering $m$-subsets $B$ of $\mathbb{Z}_n$ of the form
$$B=B_d(n,m;k_1,k_2,g,j_0)=B' \cup \bigcup_{i=1}^{c-1} (ig+H) \cup B'',$$
where $H$ is the subgroup of $\mathbb{Z}_n$ with order $d$, $g$ is an element of $\mathbb{Z}_n$, $B'$ is a proper subset of $H$ given by 
$$B'=\left\{ j \cdot \tfrac{n}{d} \; : \;  j=0,1,\dots,k_1-1 \right\},$$ and $B''$ is a proper subset of $cg+H$ of the form 
$$B''=\left\{ cg+(j_0+j) \cdot \tfrac{n}{d} \; : \;  j=0,1,\dots,k_2-1 \right\}$$ for some integer $j_0$ with $0 \leq j_0 \leq d-1$.

It turns out that our set $B$ (under some additional assumptions to be made precise in Section \ref{special sets}) has the potential to have a restricted $h$-fold sumset of size less than $u \hat{\;} (n,m,h)$ in only three cases:
\begin{itemize}
\item $h=2$, $m-1$ is not a power of 2, and $n$ is divisible by $2m-2$; 
\item $h=3$, $m=6$, and $n$ is divisible by 10; or
\item $h$ is odd, $m+2$ is divisible by $h+2$, and $n$ is divisible by $hm-h^2$.
\end{itemize}
In particular, we have the following.

\begin{prop} \label{special sets listed}
Let $B_d(n,m;k_1,k_2,g,j_0)$ be the $m$-subset of $\mathbb{Z}_n$ defined above, and let $h$ be a positive integer with $h \leq m-1$.  
\begin{itemize}
\item If $h=2$, $m-1$ is not a power of 2, $n$ is divisible by $2m-2$, and $d$ is an odd divisor of $m-1$ with $d>1$, then 
$$B_d \left(n,m;\tfrac{d+1}{2}, \tfrac{d+1}{2}, \tfrac{n}{2m-2},\tfrac{d-1}{2} \right)$$  has a restricted 2-fold sumset of size $2m-4$.
\item If $h=3$, $m=6$, and $n$ is divisible by 10, then
$$B_5 \left(n,6;4,2, \tfrac{n}{10},3 \right)$$  has a restricted 3-fold sumset of size $3m-9=9$.

\item If $h$ is odd, $m+2$ is divisible by $h+2$, and $n$ is divisible by $hm-h^2$, then
$$B_{h+2} \left(n,m;h+1,h+1, \tfrac{n}{hm-h^2},\tfrac{h+3}{2} \right)$$  has a restricted $h$-fold sumset of size $hm-h^2-1$.
\end{itemize}

\end{prop}
Our examples above demonstrate the three cases of Proposition \ref{special sets listed} in order: we have
$$C_1=B_3(12,7;2,2,1,1),$$
$$C_2=B_5(10,6;4,2,1,3),$$ and
$$C_3=B_5(15,8;4,4,1,3).$$
The proof of Proposition \ref{special sets listed} is an easy exercise; it will also follow from our proof of Theorem \ref{thm special} where we verify that, in a certain sense that we make precise, there are no other such sets. 

Let us now review some of the results that are known about the exact value of $\rho \hat{\;}(G, m, h)$.  As we have mentioned, we clearly have $\rho \hat{\;}(G, m, 1)=\rho \hat{\;}(G, m, m-1)=m$, and $\rho \hat{\;}(G, m, 1)=0$ for $h>m$; therefore, we can assume that $2 \leq h \leq m-2$.  (Since we also have $|(m-h) h\hat{\;}A| = |h\hat{\;}A|,$ we could further assume that $h \leq \left \lfloor m/2 \right \rfloor$, but we see no reason to do so.)      

Recall that, by (\ref{u-hat(p,m,h)}), for a prime $p$ we have  
$$ \rho \hat{\;}(\mathbb{Z}_p,m,h) \leq \min \{ p, hm-h^2+1\}.$$  The conjecture that equality holds here for $h=2$  
 has been known since the 1960s as the {\em Erd\H{o}s--Heilbronn Conjecture} (not mentioned in \cite{ErdHei:1964a} but in \cite{ErdGra:1980a}).  Three decades later, Dias da Silva and Hamidoune \cite{DiaHam:1994a} succeeded in proving the Erd\H{o}s--Heilbronn Conjecture and established the following.

\begin{thm} [Dias da Silva and Hamidoune; cf.~\cite{DiaHam:1994a}] \label{Dias da Silva and Hamidoune}
For a positive prime $p$ we have $$\rho\hat{\;} (\mathbb{Z}_p,m,h) = \mathrm{min} \{p, hm-h^2+1\}.$$

\end{thm}  
Theorem \ref{Dias da Silva and Hamidoune} was re-established soon after by Alon, Nathanson, and Ruzsa using their so-called polynomial-method; cf.~\cite{AloNatRuz:1995a}, \cite{AloNatRuz:1996a}, and \cite{Nat:1996a}.

For composite values of $n$, we know very little about $\rho\hat{\;} (\mathbb{Z}_n,m,h)$ in general.   One reason for this is that we do not have a version of Kneser's Theorem for restricted sumsets.

We do know a bit about the value of $\rho\hat{\;}(\mathbb{Z}_n,m,h)$ in its extremal cases.  It is easy to see that $$\rho\hat{\;}(\mathbb{Z}_n,m,h) \geq m,$$ and we can also verify that equality holds when $h=1$  or $h=m-1$, when $m \in D(n)$, or when $n$ is even, $m=4$, and $h=2$.  The fact that the inequality is strict in all other cases was proved for odd $n$ by Wang in 2008 \cite{Wan:2008a}, and for all $n$ by Girard, Griffiths, and Hamidoune in 2012 \cite{GirGriHam:2012a}.  In fact, \cite{GirGriHam:2012a} provides a complete analysis of the case $\rho\hat{\;}(G,m,h) \leq m$  for all abelian groups.

As an upper bound, we of course have  $$\rho\hat{\;}(\mathbb{Z}_n,m,h) \leq n;$$ we call the minimum value of $m$ for which equality occurs the {\em restricted $h$-critical number of $\mathbb{Z}_n$}.  It is an easy exercise to establish that the restricted 2-critical number of $\mathbb{Z}_n$ equals $\left \lfloor \tfrac{n}{2} \right \rfloor +2$; in particular, that
$$\rho \hat{\;}(\mathbb{Z}_n, \left \lfloor \tfrac{n}{2} \right \rfloor+2, 2)=n.$$  It is also not hard to see (for example, from Corollary \ref{cor h=2} in Section \ref{h=2 and h=3}) that we have 
$$\rho \hat{\;}(\mathbb{Z}_n, \left \lfloor \tfrac{n}{2} \right \rfloor+1, 2) \leq \left\{
\begin{array}{ll}
n-1 & \mbox{if} \; \exists k \in \mathbb{N}, n=2^k;\\ \\
n-2 & \mbox{otherwise.}
\end{array}
\right.$$
In 2002, Gallardo, Grekos, et al.~proved that, in fact, equality holds.

\begin{thm} [Gallardo, Grekos, et al.; cf.~\cite{GalGre:2002a}] \label{thm g G et al}
For every positive integer $n \geq 2$ we have
$$\rho \hat{\;}(\mathbb{Z}_n, \left \lfloor \tfrac{n}{2} \right \rfloor+1, 2) = \left\{
\begin{array}{ll}
n-1 & \mbox{if} \; \exists k \in \mathbb{N}, n=2^k;\\ \\
n-2 & \mbox{otherwise.}
\end{array}
\right.$$
\end{thm}

With some hesitation, we pose the following general conjecture.

\begin{conj} \label{my conj}
For all positive integers $n$ and $m$ with $3 \leq m \leq n$ we have
$$\rho\hat{\;} (\mathbb{Z}_n,m,2) = \left\{
\begin{array}{ll}
\min\{\rho(\mathbb{Z}_n,m,2), 2m-4\} & \mbox{if} \; 2|n \; \mbox{and} \; 2|m, \; \mbox{or} \\ 
& (2m-2)|n \; \mbox{and} \; \not \exists k \in \mathbb{N}, m=2^k+1; \\ \\
\min\{\rho(\mathbb{Z}_n,m,2), 2m-3\} & \mbox{otherwise.} 
\end{array}
\right.$$
\end{conj}
With Corollary \ref{cor h=2} we establish that $\rho\hat{\;} (\mathbb{Z}_n,m,2)$ cannot be more than the value given by Conjecture \ref{my conj}; and, according a computer search we performed, equality indeed holds for all $3 \leq m \leq n \leq 40$.  Note also that Conjecture \ref{my conj}, once established, would generalize Theorems \ref{Dias da Silva and Hamidoune} and \ref{thm g G et al}.    

We know less about the cases of $h \geq 3$; in particular, we do not know the reduced $h$-critical number of $\mathbb{Z}_n$ in general, in spite of some promising approaches for $h=3$ by Gallardo, Grekos, et al. in ~\cite{GalGre:2002a} and by Lev in \cite{Lev:2002a}. 

Finding $\rho \hat{\;}(G, m, h)$ for general finite abelian groups is even more difficult; as many examples demonstrate, in contrast to the case for unrestricted sumsets, the value of $\rho \hat{\;}(G, m, h)$ depends heavily on the structure of $G$ and not just its size.  To illustrate the difficulty, we point out that even $\rho \hat{\;}(\mathbb{Z}_p^2,m,2)$ is not fully known in general (cf.~the series of papers \cite{EliKer:1998a}, \cite{EliKer:2001a}, and \cite{EliKer:2001b}, by Eliahou and Kervaire).  One rare result for the exact value of $\rho \hat{\;}(G, m, h)$ is the following extension of the Erd\H{o}s--Heilbronn problem.

\begin{thm} [K\'arolyi; cf.~\cite{Kar:2004a}, \cite{Kar:2003a}] 
Suppose that the smallest prime divisor of $|G|$ is at least $2m-3$.  Then $$\rho \hat{\;}(G, m, 2)=2m-3.$$
\end{thm}   

We should also mention the following conjecture of Plagne.  (See also similar conjectures of Lev in \cite{Lev:2000a} and \cite{Lev:2005a}.)

\begin{conj} [Plagne; cf.~\cite{Pla:2006b}] \label{Plagne's conj}
For every abelian group $G$ and for every $m \leq |G|$ we have
$$\rho \hat{\;}(G, m, 2) \geq \rho (G, m, 2)-2.$$

\end{conj}
We should point out that Conjecture \ref{my conj} implies Conjecture \ref{Plagne's conj} for cyclic groups.  Indeed, we have $$\rho (\mathbb{Z}_n,m,2) \leq f_1=2m-1;$$ furthermore, when $n$ is even, we have $$\rho (\mathbb{Z}_n,m,2) \leq f_2=2m-2,$$ and when $n$ is divisible by $2m-2$, we find that $$\rho (\mathbb{Z}_n,m,2) \leq f_{2m-2}=2m-2.$$

However, a similar bound for $h \geq 3$ will definitely not hold; as the following proposition demonstrates, the difference of $\rho (G, m, h)$ and $\rho \hat{\;}(G, m, h)$ may be arbitrarily large. 

\begin{prop} \label{rho-hat very small}
For every $h \geq 3$ and for every positive real number $C$, one can find a finite abelian group $G$ and a positive integer $m$ so that
$$\rho \hat{\;}(G, m, h) < \rho (G, m, h)-C.$$
\end{prop}

{\em Proof}.  Choose a prime $p>h$ and a positive integer $t >2$ so that $$(h-2)p^{t-1} \geq C.$$
One can readily verify that for $n=p^t$ and $m=p^{t-1}+1$ we have 
$$\rho (\mathbb{Z}_n, m, h)=u(n,m,h)=f_1=hp^{t-1}+1$$ and
$$\rho \hat{\;}(\mathbb{Z}_n, m, h) \leq u \hat{\;} (n,m,h)=f \hat{_{p^{t-1}}}=2p^{t-1},$$ from which our claim follows.
$\Box$

Before ending this section, we exhibit some counterexamples for conjectures that have appeared in the literature.  

In their interesting paper \cite{HamLlaSer:2000a}, Hamidoune, Llad\'o, and Serra conjectured (Conjecture 3.5) that if $A$ is a finite generating $m$-subset of an abelian group $G$ of order $n$ with $0 \in A$ and $m \geq 6$, then 
$$|2 \hat{\;} A| \geq \min \{n-1,3m/2\}.$$
They mention that $m \geq 6$ is necessary, since if $A$ is an arithmetic progression of length $m<6$, then $|2 \hat{\;} A|=2m-3 < 3m/2$.  However, one can see that there are $7$-subsets of $\mathbb{Z}_{12}$ contradicting this claim; for example, for the set $C_1$ mentioned above, we have $$|2 \hat{\;} C_1|=10 < \min \{12-1,3 \cdot 7 /2\}.$$ 

In fact, one can find arbitrarily large counterexamples, as follows.  Suppose that $G$ is of the form 
$$G=\mathbb{Z}_2^k \oplus \mathbb{Z}_d$$ with $k \geq 1$ and $d \geq 3$ odd.  Consider the subset $$A=\mathbb{Z}_2^k \oplus \{0,1\}$$ of $G$.  Observe that $A$ generates $G$, and 
$$2 \hat{\;} A=(\mathbb{Z}_2^k \oplus \{0,1,2\}) \setminus ({\bf 0} \oplus \{0,2\}),$$ where ${\bf 0}$ is the identity element of $\mathbb{Z}_2^k$.
Then $m=|A|=2 \cdot 2^k$ and $|2 \hat{\;} A|=3 \cdot 2^k-2$, contradicting the conjecture.

Plagne, in his powerful paper \cite{Pla:2006b}, made the following conjecture (as part of Conjecture 9):  For every abelian group $G$ of odd order $n$, there is a constant $c_G$, so that $c_G=o(n)$ as $n$ tends to infinity, and for every $m$ with $c_G \leq m \leq n$ that is ``exceptional'' for $G$, we have $$\rho \hat{\;}(G, m, 2) = \min \{\rho (G, m, 2), 2m-2\}.$$  Plagne defines $m$ to be exceptional for $G$ if for every finite nontrivial minimal subgroup $M$ of $G$ one has $m \equiv 1$ mod $|M|$.  

For a counterexample, let $p$ be an odd prime, $t \geq 2$, $n=p^t$, $G=\mathbb{Z}_n$, and $m=p^{t-1}+1$ (cf.~the proof of Proposition \ref{rho-hat very small}).  Then $m$ is exceptional for $G$; furthermore, $$\lim_{n \rightarrow \infty} m/n = 1/p>0.$$  It can be easily verified that $$\rho (G, m, 2)= u(n,m,2)=f_1=2m-1$$ but $$\rho \hat{\;}(G, m, 2) \leq u\hat{\;}(n,m,2) = f\hat{_1}=2m-3.$$ 
Therefore, $$\rho \hat{\;}(G, m, 2) \neq  \min \{\rho (G, m, 2), 2m-2\}.$$

The pursuit of finding the value of $\rho \hat{\;}(G, m, h)$ remains challenging and exciting.

\section{The set $A_d(n,m)$}  \label{regular sets}

In this section we compute $|h A_d(n,m)|$ and $|h\hat{\;}A_d(n,m)|$; in particular, we prove (\ref{|hA_d(n,m)|}) and Theorem \ref{sizeofh^A_d}.  But first, a brief justification for why the set $A_d(n,m)$ is of interest.

How can one find $m$-subsets $A$ in a group $G$ that have small $h$-fold sumsets $hA$?  Two ideas come to mind.  First, observe that if $A$ is a subset of a subgroup $H$ of $G$, then $hA$ will be a subset of $H$ as well; a bit more generally, if $A$ is a subset of the coset $g+H$ of $H$ for some $g \in G$, then $hA$ will be a subset of the coset $hg+H$.  Therefore, we have $|hA| \leq d$ for every divisor $d$ of $n$ that is not less than $m$.

The second idea is based on the observation that, when $A$ is an arithmetic progression
$$A=\{a,a+g,a+2g,\dots,a+(m-1)g\}$$ for some $a, g \in G$, then many of the $h$-fold sums coincide; in particular, 
$$hA=\{ha,ha+g,ha+2g,\dots,ha+(hm-h)g\},$$  thus, consequently, we have $|hA| \leq hm-h+1$.

As a combination of these two ideas, we choose a subgroup $H$ of $G$, and then select an $m$-subset $A$ of $G$ so that its elements are in as few cosets of $H$ as possible; furthermore, we want these cosets to form an arithmetic progression.   

More explicitly, we let $G=\mathbb{Z}_n$, fix a divisor $d$ of $n$, and consider the (unique) subgroup $H$ of order $d$ of $\mathbb{Z}_n$,  namely, $$H=\left \{j \cdot \tfrac{n}{d}  \; : \;  j=0,1,2,\dots,d-1 \right\}.$$      
With $|A|=m$ and $|H|=d$, the number of cosets of $H$ that we need is $\left \lceil \tfrac{m}{d} \right \rceil$.  We let $k$ denote the positive remainder of $m$ mod $d$; that is, $k=m-dc$, where
$c=\left \lceil \tfrac{m}{d} \right \rceil-1.$  Now we choose our cosets to be
$i +H$ with $i=0,1,2,\dots, c$, and set
$$A_d=A_d(n,m)=\bigcup_{i=0}^{c-1} (i+H) \cup \left\{c + j \cdot \tfrac{n}{d}  \; : \;  j=0,1,2,\dots,k-1 \right\}.$$ 

Then $m \leq n$ assures that $c < n/d$, thus $A_d$ has $m$ distinct elements.  Note that, when $m \leq d$ (that is, $c=0$), then $A_d$ lies entirely within a single coset and forms an arithmetic progression.  We should also point out that here we chose the cosets represented by $a=0$ and $g=1$; the assumption $a=0$ we can make without any loss of generality since $|h(a+A)|$ does not depend on $a$, and we also chose $g=1$ as the general case does not yield any benefits here (in contrast to Section \ref{special sets}).

We see that
$$hA_d(n,m)=\bigcup_{i=0}^{hc-1} (i+H) \cup \left\{hc + j \cdot \tfrac{n}{d}  \; : \;  j=0,1,2,\dots,h(k-1) \right\},$$
and thus
$$
|hA_d(n,m)| = \min \{ n , hcd + \min \{d, h(k-1)+1 \}; $$
since $hcd+d=f_d$ and $hcd+h(k-1)+1=hm-h+1$,
we get (\ref{|hA_d(n,m)|}), from which Theorem \ref{value of u} follows for cyclic groups as explained in the Introduction.

We now turn to the computation of $
|h \hat{\;} A_d(n,m)|$ and the proof of Theorem \ref{sizeofh^A_d}.
We first state and prove the following lemma.

\begin{lem} \label{lemma any j}
Suppose that $d$ and $t$ are positive integers with $t \leq d-1$, and let $j \in \mathbb{Z}_d$.  Then there is a $t$-subset $J=\{j_1,\dots,j_t\}$ of $\mathbb{Z}_d$ for which $$j_1+j_2+\cdots+j_t=j.$$  

\end{lem}

Note that the restriction of $t \leq d-1$ is necessary: for $t=d$ there is exactly one $j \in \mathbb{Z}_d$ for which such a set exists, and there are no such sets when $t>d$.

{\em Proof of Lemma \ref{lemma any j}.}  Let $j_0$ be the nonnegative remainder of the integer $$j-\tfrac{t^2-t}{2}$$ mod $d$.  We can then easily check that the set $J$ defined (for example) as  
$$
J = \left\{ 
\begin{array}{ll}
\{0,1,\dots,t-1,t\} \setminus \{t-j_0\}  & \mbox{if $0 \leq j_0 \leq t-1$,}\\ \\
\{1,\dots,t-1\} \cup \{j_0\}  & \mbox{if $t \leq j_0 \leq d-1$}\\ 
\end{array}\right.
$$ satisfies our requirements.
$\Box$

{\em Proof of Theorem \ref{sizeofh^A_d}.}  Let us recall our notations and introduce some new ones.  We write
$$m=dc+k \; \mbox{with} \; c=\left \lceil \tfrac{m}{d}  \right \rceil -1$$
and
$$h=dq+r \; \mbox{with} \; q=\left \lceil \tfrac{h}{d}  \right \rceil -1;$$
note that $ c \geq 0$, $q \geq 0$, $1 \leq k \leq d$, and $1 \leq r \leq d$.  Recall also that we have set $$A_d=A_d(n,m)=\bigcup_{i=0}^{c-1} (i+H) \cup \left\{c + j \cdot \tfrac{n}{d}  \; : \;  j=0,1,2,\dots,k-1 \right\}$$ where $H$ is the order $d$ subgroup of $\mathbb{Z}_n$.
Note that every element of $h \hat{\;} A_d$ is of the form 
$$(i_1+i_2+\cdots+i_h)+(j_1+j_2+\cdots+j_h) \cdot \tfrac{n}{d}$$
with $i_1,\dots,i_h \in \{0,1,\dots,c\}$ and $j_1,\dots,j_h \in \{0,1,\dots,d-1\}$, with the added conditions that when any of the $i$-indices equals $c$, the corresponding $j$-index is at most $k-1$, and that when two $i$-indices are equal, the corresponding $j$-indices are distinct.

Clearly, the least value  of $i_1+\dots+i_h$ is 
\begin{eqnarray} \label{i min}
i_{\mathrm{min}}=d(0+1+\cdots+(q-1))+rq=q (h+r-d)/2.
\end{eqnarray} 

To compute the largest value $i_{\mathrm{max}}$ of $i_1+\dots+i_h$, we consider four cases depending on whether $r >k$ or not and whether $q=0$ or not.  

First, when $q=0$ and $r>k$, then $r=h$ and $1 \leq h-k \leq d$, so it is easy to see that
$$i_{\mathrm{max}}=kc+(h-k)(c-1)=hc-h+k.$$
Next, when  $q=0$ and $r \leq k$, then $r=h \leq k$, so we have
$$i_{\mathrm{max}}=hc=hc-h+r.$$
In the case when $r >k$ and $q \geq 1$, we write $h$ as $h=k+dq+(r-k)$; thus
\begin{eqnarray*}
i_{\mathrm{max}} & = & k c+[(c-1)+(c-2)+\cdots+(c-q)]d+(r-k)(c-q-1) \\
& = & hc-h+kq-q (h+r-d)/2+k.
\end{eqnarray*} 
Finally, in the case when $q \geq 1$ and $r \leq k$, then $1 \leq d-k+r  \leq d$; we write $h$ as $h=k+d(q-1)+(d-k+r)$ and thus (using our result from the previous case)
\begin{eqnarray*}
i_{\mathrm{max}} & = & k c+[(c-1)+(c-2)+\cdots+(c-q+1)]d+(d+r-k)(c-q) \\
& = & hc-h+kq-q (h+r-d)/2+r.
\end{eqnarray*} 
All four cases can be summarized by the formula
\begin{eqnarray} \label{i max}
i_{\mathrm{max}}=hc-h+kq-q (h+r-d)/2+ \min\{r,k\}.
\end{eqnarray}

Obviously, $i=i_1+i_2+\cdots+i_h$ can assume the value of any integer between these two bounds, and thus $h \hat{\;} A_d$ lies in exactly $$\min \left\{n/d,i_{\mathrm{max}}-i_{\mathrm{min}}+1 \right\}$$ cosets of $H$.  This immediately yields the upper bound 
$$|h \hat{\;} A_d| \leq   \left\{n,\left(i_{\mathrm{max}}-i_{\mathrm{min}}+1 \right) d \right\} ,$$
where, by a simple calculation from (\ref{i min}) and (\ref{i max}) above, 
\begin{eqnarray} \label{imax-imin+1}
\left(i_{\mathrm{max}}-i_{\mathrm{min}}+1 \right) d = hm-h^2-r(k-r)+d \min\{0,k-r\}+d.
\end{eqnarray}

We will separate the rest of the proof into several cases.  We can easily check (by considering the cases of $d=1$ and $d \geq 2$ separately) that our formula holds for $h=1$, so below we assume that $2 \leq h \leq m-1$.

{\bf Claim 1}: If $h \leq k$ and $h<d$, then $|h \hat{\;} A_d| = \min\{n, f_d, hm-h^2+1\}$.

{\em Proof of Claim 1}:  Note that the assumptions give $i_{\mathrm{min}}=0$, $i_{\mathrm{max}}=hc$; using Lemma \ref{lemma any j} above, we have
$$h \hat{\;} A_d =\bigcup_{i=0}^{hc-1} (i+H) \cup \left\{hc+j \cdot \tfrac{n}{d}  \; : \;  j=\tfrac{h(h-1)}{2}, \dots, h(k-1)-\tfrac{h(h-1)}{2} \right\}.$$
Therefore, we see that
\begin{eqnarray*}
|h \hat{\;} A_d| & =& \min \{ n , hcd + \min \{d, h(k-1)-h(h-1)+1 \} \\ 
& =& \min \{ n , hcd+d , hm-h^2+1 \} \\
& =& \min \{ n , f_d , hm-h^2+1 \},
\end{eqnarray*}
as claimed.

{\bf Claim 2}: If $h=k=d$, then $|h \hat{\;} A_d| = \min\{n, hm-h^2-h+2\}$.

{\em Proof of Claim 2}:  Our assumptions yield $i_{\mathrm{min}}=0$, $i_{\mathrm{max}}=hc=dc$; we note that $h<m=dc+k=dc+h$, so $c \geq 1$.  In the case of $h=k=d$, Lemma \ref{lemma any j} cannot be used for the coset $i+H$ when $i=0$; we now have  
$$h \hat{\;} A_d =\left\{\tfrac{d(d-1)}{2} \cdot \tfrac{n}{d} \right\} \cup \bigcup_{i=1}^{dc-1} (i+H) \cup \left\{dc+ \tfrac{d(d-1)}{2} \cdot \tfrac{n}{d}  \right\}.$$

If we have $$dc-1 \geq \tfrac{n}{d},$$ then clearly $|h \hat{\;} A_d|=n$, but we also have
$$hm-h^2-h+2=h(dc+h)-h^2-h+2=h(dc-1)+2 \geq n+2>n,$$
so 
$$|h \hat{\;} A_d| = \min\{n, hm-h^2-h+2\},$$ as claimed.

Assume now that 
 $$dc-1 \leq \tfrac{n}{d}-2.$$
In this case
$$|h \hat{\;} A_d| =1+(dc-1)d+1 = h(cd+h)-h^2-h+2=hm-h^2-h+2;$$ furthermore, 
$$1+(dc-1)d+1 \leq n-2d+2  \leq n,$$  thus
$$|h \hat{\;} A_d| = \min\{n, hm-h^2-h+2\},$$ as claimed.

This leaves us with the case of $$dc-1 = \tfrac{n}{d}-1,$$  when we have
\begin{eqnarray*}
h \hat{\;} A_d &=& \left\{\tfrac{d(d-1)}{2} \cdot \tfrac{n}{d} \right\} \cup \bigcup_{i=1}^{dc-1} (i+H) \cup \left\{dc+ \tfrac{d(d-1)}{2} \cdot \tfrac{n}{d}  \right\} \\ \\
&=& \left\{\tfrac{d(d-1)}{2} \cdot \tfrac{n}{d} \right\} \cup \bigcup_{i=1}^{\tfrac{n}{d}-1} (i+H) \cup \left\{\tfrac{n}{d}+ \tfrac{d(d-1)}{2} \cdot \tfrac{n}{d}  \right\} \\ \\
&=& \left\{j \cdot \tfrac{n}{d} \; : \;  j= \tfrac{d(d-1)}{2}, \tfrac{d(d-1)}{2}+1 \right\} \cup \bigcup_{i=1}^{\tfrac{n}{d}-1} (i+H).
\end{eqnarray*}
But $d=h \geq 2$, so, as above, we find that
$$|h \hat{\;} A_d|=(dc-1)d+2 =hm-h^2-h+2.$$
Furthermore, $2+\left(\tfrac{n}{d}-1 \right) \cdot d=n-(d-2) \leq n,$ so again we have
$$|h \hat{\;} A_d| = \min\{n, hm-h^2-h+2\},$$ completing the proof of Claim 2.

{\bf Claim 3}: If $h>k$, $r \neq d$, and $r \neq k$, then $$|h \hat{\;} A_d| = \min\{n, hm-h^2-r(k-r)+d \min\{0,k-r\}+d\}.$$

{\em Proof of Claim 3}:  Observe that, by Lemma \ref{lemma any j}, the three conditions imply that 
$$h \hat{\;} A_d=\bigcup_{i=i_{\mathrm{min}}}^{i_{\mathrm{max}}} (i+H).$$  Our result now follows from (\ref{imax-imin+1}).

{\bf Claim 4}: If $h>k$, $r=d$, and $r \neq k$, then $$|h \hat{\;} A_d| = \min\{n, hm-h^2+1\}.$$

{\em Proof of Claim 4}:  This time we have 
$$h \hat{\;} A_d=\{ x_{\mathrm{min}} \} \cup \bigcup_{i=i_{\mathrm{min}}+1}^{i_{\mathrm{max}}} (i+H)$$
where $x_{\mathrm{min}}$ equals the sum of the $h$ elements of the set
$$\bigcup_{i=0}^{q} (i+H).$$  
Therefore, 
$$|h \hat{\;} A_d|=\min\{n, (i_{\mathrm{max}}-i_{\mathrm{min}})\cdot d +1 \}.$$
Our claim then follows, since we now have $k<r=d$ and thus from (\ref{imax-imin+1}) we get
$$(i_{\mathrm{max}}-i_{\mathrm{min}})\cdot d +1=(i_{\mathrm{max}}-i_{\mathrm{min}}+1)\cdot d -(d-1)=hm-h^2+1.$$

{\bf Claim 5}: If $h>k$, $r \neq d$, and $r = k$, then $$|h \hat{\;} A_d| = \min\{n, hm-h^2+1\}.$$

{\em Proof of Claim 5}:  Our conditions imply that
$$h \hat{\;} A_d= \bigcup_{i=i_{\mathrm{min}}}^{i_{\mathrm{max}}-1} (i+H) \cup \{ x_{\mathrm{max}} \}$$
where $x_{\mathrm{max}}$ equals the sum of the $h$ elements of the set
$$\bigcup_{i=c-q}^{c-1} (i+H) \cup \left\{c + j \cdot \tfrac{n}{d}  \; : \;  j=0,1,2,\dots,k-1 \right\}.$$  Our claim then follows as in Claim 4.

{\bf Claim 6}: If $h>k$, $r =k= d$, then 
$$|h \hat{\;} A_d| = \min\{n, hm-h^2-d+2\}.$$

{\em Proof of Claim 6}:  This time we get
$$h \hat{\;} A_d=\{ x_{\mathrm{min}} \} \cup \bigcup_{i=i_{\mathrm{min}}+1}^{i_{\mathrm{max}}-1} (i+H) \cup \{ x_{\mathrm{max}} \}$$
where $x_{\mathrm{min}}$ and $x_{\mathrm{max}}$ were defined in the proofs of Claims 4 and 5, respectively.

With $r=k=d$, we see from (\ref{imax-imin+1}) above that
$$i_{\mathrm{max}}-i_{\mathrm{min}}-1 = \tfrac{hm-h^2-d}{d}.$$

We consider three subcases.  If $m$ is more than $n/h+h$, then $i_{\mathrm{max}}-i_{\mathrm{min}}-1$ is more than $n/d-1$ so, since it is an integer, it is at least $n/d$.  Therefore, $|h \hat{\;} A_d| =n$; furthermore, 
$$\min\{n, hm-h^2-d+2 \}= \min \{n, \left(i_{\mathrm{max}}-i_{\mathrm{min}}-1 \right)d+2 \} \geq \min\{n, n+2   \}=n.$$

If $m$ is less than $n/h+h$, then $i_{\mathrm{max}}-i_{\mathrm{min}}-1$ is less than $n/d-1$ and thus at most $n/d-2$.  Therefore, 
$$|h \hat{\;} A_d| =  hm-h^2-d+2;$$
furthermore, $hm-h^2-d+2$ is less than $n-d+2$ and thus at most $n-d+1$, which is never more than $n$, and thus again we have
$$|h \hat{\;} A_d| =  \min \{n, hm-h^2-d+2\}.$$
 
This leaves us with the case of $m=\tfrac{n}{h}+h$, in which case our claim becomes
$$|h \hat{\;} A_d| = \min\{n,n-d+2\}.$$

We now have
\begin{eqnarray*}
h \hat{\;} A_d & =& \{ x_{\mathrm{min}} \} \cup \bigcup_{i=i_{\mathrm{min}}+1}^{i_{\mathrm{max}}-1} (i+H) \cup \{ x_{\mathrm{max}} \} \\
& =& \bigcup_{i=i_{\mathrm{min}}+1}^{i_{\mathrm{min}}+\tfrac{n}{d}-1} (i+H) \cup \{x_{\mathrm{min}},  x_{\mathrm{max}} \}.
\end{eqnarray*}

A simple computation shows that, denoting the sum of the elements in a subset $S$ of $\mathbb{Z}_n$ by $\sum S$, we have
$$x_{\mathrm{min}} = \sum \bigcup_{i=0}^{q} (i+H) = \tfrac{dq(q+1)}{2} + \tfrac{d(d-1)(q+1)}{2} \cdot \tfrac{n}{d}$$
and
$$x_{\mathrm{max}} =\sum \bigcup_{i=c-q}^{c} (i+H)  = cd(q+1) - \tfrac{dq(q+1)}{2} + \tfrac{d(d-1)(q+1)}{2} \cdot \tfrac{n}{d}.$$
But
$$cd=m-d=\tfrac{n}{h}+h-d=\tfrac{n}{d(q+1)}+dq,$$
thus
$$x_{\mathrm{max}} = \tfrac{dq(q+1)}{2} + \left( \tfrac{d(d-1)(q+1)}{2} +1 \right) \cdot \tfrac{n}{d},$$
showing that $x_{\mathrm{min}} = x_{\mathrm{max}}$ if, and only if, $d=1$.
Therefore, when $d \geq 2$, we get
$$|h \hat{\;} A_d|=\left(\tfrac{n}{d}-1 \right)d+2=n-d+2,$$ and when $d=1$ we get $$|h \hat{\;} A_d|=\left(\tfrac{n}{d}-1 \right)d+1=n.$$
This completes the proof of Claim 6 and, with that, Theorem \ref{sizeofh^A_d} is established.
$\Box$

\section{Special sets}  \label{special sets}

We here present a variation on the construction of Section \ref{regular sets} that, under specific conditions, yields smaller sumsets.  The subsets we present here are very similar to $A_d(n,m)$: we still pack our $m$  elements into the $\lceil m/d \rceil $ cosets of the subgroup $H$ of $\mathbb{Z}_n$ of order $d$, these cosets are still forming an arithmetic progression, but this time we have the common difference $g$ that is not necessarily 1, and we allow not one but two of the cosets to not be entirely in our set.  

To do this, we write $m$ in the form
$$m=k_1+(c-1)d+k_2$$ for some positive integers $c$, $k_1$, and $k_2$; we assume that $k_1<d$, $k_2<d$, but $k_1+k_2>d$, and thus, as before, $$c=\lceil m/d \rceil -1.$$  We then are considering $m$-subsets $B$ of $\mathbb{Z}_n$ of the form
$$B=B' \cup \bigcup_{i=1}^{c-1} (ig+H) \cup B'',$$
where $H$ is the subgroup of $\mathbb{Z}_n$ with order $d$, $g$ is an element of $\mathbb{Z}_n$, $B'$ is a proper subset of $H$ given by 
$$B'=\left\{ j \cdot \tfrac{n}{d} \; : \;  j=0,1,\dots,k_1-1 \right\},$$ and $B''$ is a proper subset of $cg+H$ of the form 
$$B''=\left\{ cg+(j_0+j) \cdot \tfrac{n}{d} \; : \;  j=0,1,\dots,k_2-1 \right\}$$ for some integer $j_0$ with $0 \leq j_0 \leq d-1$.

For a nonnegative integer $i$, we let $$h \hat{_i} B=h \hat{\;} B \cap (ig+H),$$ with which we can write 
$$h \hat{\;} B = \bigcup_{i=i_{\min}}^{i_{\max}} h \hat{_i} B \subseteq \bigcup_{i=i_{\min}}^{i_{\max}} (ig+H),$$ where $i_{\min}$ and $i_{\max}$ denote the value of $$i_1+i_2+\cdots+i_h$$ for the ``first'' and ``last'' $h$ elements of $B$, respectively.  Note that $c \geq 1$ implies that $i_{\max} > i_{\min}$.  We will select $g \in \mathbb{Z}_n$ so that
\begin{eqnarray} \label{g def}
i_{\max} \cdot g= i_{\min} \cdot g+ \tfrac{n}{d};
\end{eqnarray}
with this choice we have $$i_{\max} \cdot g+H=i_{\min} \cdot g+H.$$  Therefore, $h \hat{_{i_{\max}}} B$ and $h \hat{_{i_{\min}}} B$ are subsets of the same coset of $H$; we are interested in the situation when one of them is a subset of the other, in which case we call $B$ {\em special}. 

In particular, we examine special sets in two cases: when both $k_1$ and $k_2$ are greater than or equal to $h$, and when exactly one of them is less than $h$.  The extent to which such special sets yield an improvement over the upper bound $u \hat{\;} (n,m,h)$ in these cases, and the conditions under which they occur, are given by the following theorem. 

\begin{thm} \label{thm special}
Suppose that the set $B$ (as defined above) is special and that $h \leq \max\{k_1,k_2\}$.   Then
$$|h \hat{\;} B| \geq u \hat{\;} (n,m,h) -1;$$
furthermore, if equality holds, then 
\begin{itemize}
\item $h=2$, $m-1$ is not a power of 2, and $n$ is divisible by $2m-2$; 
\item $h=3$, $m=6$, and $n$ is divisible by 10; or
\item $h$ is odd, $m+2$ is divisible by $h+2$, and $n$ is divisible by $hm-h^2$.
\end{itemize}

\end{thm}

Our proof below will also verify Proposition \ref{special sets listed}.

{\em Proof.}  Without loss of generality, we assume that $k_2 \leq k_1$, in which case $B$ is special when  
\begin{eqnarray}
h \hat{_{i_{\max}}} B \subseteq h \hat{_{i_{\min}}} B. 
\end{eqnarray}
We will consider to cases: when $h \leq k_2 \leq k_1$ and when $k_2<h \leq k_1$.

{\bf Case 1}: Suppose first that $B$ is a special set and that $h \leq k_2 \leq k_1$, in which case
$i_{\min}=0$ and  $i_{\max}=hc$; furthermore, to satisfy (\ref{g def}), $n$ must be divisible by $hcd$, in which case we set $$g=\tfrac{n}{hcd}.$$ 
We thus see that
$$h \hat{_{i_{\min}}} B=h \hat{\;} B'=\left\{j \cdot \tfrac{n}{d}  \; : \;  j=\tfrac{h(h-1)}{2}, \dots, h(k_1-1)-\tfrac{h(h-1)}{2} \right\},$$ and 
\begin{eqnarray*}
h \hat{_{i_{\max}}} B=h \hat{\;} B'' & = & \left\{hcg+(hj_0+j) \cdot \tfrac{n}{d}  \; : \;  j=\tfrac{h(h-1)}{2}, \dots, h(k_2-1)-\tfrac{h(h-1)}{2} \right\} \\ \\
& = & \left\{(1+hj_0+j) \cdot \tfrac{n}{d}  \; : \;  j=\tfrac{h(h-1)}{2}, \dots, h(k_2-1)-\tfrac{h(h-1)}{2} \right\}.
\end{eqnarray*}

We start by proving that 
\begin{eqnarray}  \label{c=1 or c>1}
h \hat{\;} B = h \hat{\;} B' \cup \bigcup_{i=1}^{hc-1} (ig+H) \cup h \hat{\;} B''.
\end{eqnarray}

When $c \geq 2$, then (\ref{c=1 or c>1}) follows immediately from Lemma \ref{lemma any j}, so assume that $c=1$. In this case $B=B' \cup B''$, and
\begin{eqnarray*}
h \hat{\;} B & = & \bigcup_{i=0}^h (h-i) \hat{\;} B' + i \hat{\;} B'' \\ 
& = & \bigcup_{i=0}^h \left\{ i  g+ j \cdot \tfrac{n}{d} \; : \;  j=j_{\min}(i), j_{\min}(i)+1, \ldots, j_{\max}(i) \right\},
\end{eqnarray*}
where
$$j_{\min} (i)= \tfrac{(h-i)(h-i-1)}{2} + i j_0+ \tfrac{i(i-1)}{2},$$ and
$$j_{\max} (i)= (h-i)(k_1-1)-\tfrac{(h-i)(h-i-1)}{2} + i j_0+ i(k_2-1)-\tfrac{i(i-1)}{2}.$$
We see that, for $1 \leq i \leq h-1$, we have
\begin{eqnarray*}
j_{\max} (i)-j_{\min} (i) +1& = & (h-i)(k_1-1)+i(k_2-1)-(h-i)(h-i-1)-i(i-1)+1 \\
& \geq & (h-i)(k_1-1)+i(k_2-1)-(h-i)(h-i-1)-i(i-1) +1+d- (k_1+k_2-1) \\
& = & (h-i-1)(k_1-h+i)+(i-1)(k_2-i)-(h-2)+d\\ 
& \geq & (h-i-1)i+(i-1)(h-i)-(h-2)+d\\ 
& = & 2(i-1)(h-i-1)+d \\
& \geq & d. 
\end{eqnarray*}
Therefore,
$$h \hat{\;} B = h \hat{\;} B' \cup \bigcup_{i=1}^{h-1} (ig+H) \cup h \hat{\;} B'',$$
as claimed.

Recall that $B$ is special and thus
\begin{eqnarray} \label{special--1} 
h \hat{\;} B'' \subseteq h \hat{\;} B',
\end{eqnarray}  
so from (\ref{c=1 or c>1}) we see that
\begin{eqnarray*} 
h \hat{\;} B = h \hat{\;} B' \cup \bigcup_{i=1}^{hc-1} (ig+H),
\end{eqnarray*}
and, therefore,
\begin{eqnarray} \label{size of h^B}
|h \hat{\;} B|&=&\min\{d, hk_1-h^2+1\} + (hc-1) \cdot d.
\end{eqnarray}

Assume now that $$|h \hat{\;} B| < u \hat{\;} (n,m,h).$$  
Note that $d_0=hcd \in D(n)$ and $$m=k_1+(c-1)d+k_2 < (c+1)d \leq hcd=d_0,$$ and thus
$$u \hat{\;} (n,m,h) \leq u (n,m,h) \leq f_{d_0} = \left( h \left \lceil \tfrac{m}{d_0}   \right \rceil -h+1 \right)d_0=hcd.$$  Therefore, by (\ref{size of h^B}) we have
$$hk_1-h^2+1 < d$$ and thus
\begin{eqnarray} \label{|h^B| < hcd}
|h \hat{\;} B| &=& hk_1-h^2+1 + (hc-1) \cdot d \\
&=& hm-h^2+1-(hk_2-hd+d).  \label{|h^B| < hm-h^2+1}
\end{eqnarray}

Now if we were to have $$|h \hat{\;} B| \leq u \hat{\;} (n,m,h)-2,$$ then
$$|h \hat{\;} B| \leq \min\{hcd, hm-h^2+1\}-2,$$
so by (\ref{|h^B| < hcd}) and (\ref{|h^B| < hm-h^2+1}) we have 
$$hd-d+2 \leq hk_2 \leq hk_1 \leq h^2+d-3.$$
But, since $d>h$,  this is impossible: the left hand side is more than the right hand side when $d \geq h+2$, and for $d=h+1$, there are no multiples of $h$ between the two sides.  Therefore, $$|h \hat{\;} B| \geq u \hat{\;} (n,m,h)-1,$$ as claimed.

Assume now that $$|h \hat{\;} B| = u \hat{\;} (n,m,h)-1;$$ from which we get 
\begin{eqnarray} \label{k_1 and k_2}
hd-d+1 \leq hk_2 \leq hk_1 \leq h^2+d-2.
\end{eqnarray}

We see that (\ref{k_1 and k_2}) can only hold if $k_1=k_2$.  Indeed, the difference of the right hand side and the left hand side is at most $h-1$, since $2 \leq h \leq d-1$ implies that
$$(h^2+d-2)-(hd-d+1)-(h-1)=(h-2)(h-(d-1)) \leq 0.$$
But $k_1=k_2$ means that $h \hat{\;} B'=h \hat{\;} B'',$ which can only happen if 
\begin{eqnarray} \label{j_0 mod}
1+hj_0 \equiv 0 \; \; \mbox{mod} \; d.
\end{eqnarray} As an immediate consequence, we see that $\gcd(d,h)=1$.

Let us return to (\ref{k_1 and k_2}), which we can now write as
\begin{eqnarray} \label{k_1 = k_2}
hd-d+1 \leq hk_2 = hk_1 \leq h^2+d-2.
\end{eqnarray}

If $3 \leq h \leq d-3$, then the left hand side is more than the right hand side, unless $h=3$ and $d=6$, which leads to no solutions for $k_1$ and $k_2$.

If $h=2$, then (\ref{k_1 = k_2}) becomes $$d+1 \leq 2k_2 = 2k_1 \leq d+2;$$ since $d$ must be odd, we get $k_1=k_2=(d+1)/2$ and $m=cd+1$, so $n$ is divisible by $hcd=2m-2$.   Furthermore, $d$ is odd and $d >h=2$, so $2m-2$, and thus $m-1$, is not a power of 2.  From (\ref{j_0 mod}) we get $j_0=(d-1)/2$, and so we arrived at the set $$B_d \left(n,m;\tfrac{d+1}{2}, \tfrac{d+1}{2}, \tfrac{n}{2m-2},\tfrac{d-1}{2} \right)$$ exhibited in the first part of Proposition \ref{special sets listed}.  This set can be written explicitly as 
$$\left\{ j \cdot \tfrac{n}{d} \; : \;  j=0,1,\dots,\tfrac{d-1}{2} \right\} \cup \bigcup_{i=1}^{c-1} \left( \tfrac{i \cdot n}{2m-2} + H \right) \cup \left\{\tfrac{c \cdot n}{2m-2} + \left( \tfrac{d-1}{2} + j \right) \cdot \tfrac{n}{d} \; : \;  j=0,1,\dots, \tfrac{d-1}{2}\right\},$$ and it has a 2-fold restricted sumset of size $2m-4$.

If $h=d-2$, then the only solution to (\ref{k_1 = k_2}) is $k_1=k_2=d-1$, in which case $m=(c+1)d-2=(c+1)(h+2)-2$, so $m+2$ is divisible by $h+2$, and $n$ is divisible by $hcd=hm-h^2$.   Furthermore, note that $h=d-2$ and $\gcd(d,h)=1$ imply that $d$ and $h$ are odd.  We have $j_0=(d+1)/2$ and thus get the third part of Proposition \ref{special sets listed}, namely, $$B_{h+2} \left(n,m;h+1,h+1, \tfrac{n}{hm-h^2},\tfrac{h+3}{2} \right),$$ or  
$$\left\{ j \cdot \tfrac{n}{d} \; : \;  j=0,1,\dots,d-2 \right\} \cup \bigcup_{i=1}^{c-1} \left( \tfrac{i \cdot n}{hm-h^2} + H \right) \cup \left\{\tfrac{c \cdot n}{hm-h^2} + \left( \tfrac{d+1}{2} + j \right) \cdot \tfrac{n}{d} \; : \;  j=0,1,\dots,d-2\right\}.$$ This set has $h$-fold restricted sumset of size $hm-h^2-1$.   

Finally, if $h=d-1$, then the only solution to (\ref{k_1 = k_2}) is $k_1=k_2=d-1$, in which case $m=(c+1)d-2=(c+1)(h+1)-2$, so $m+2$ is divisible by $h+1$, and $n$ is divisible by $hcd=hm-h^2+h$.  In this case, we can assume that $h \geq 3$, since for $h=1$ we always have $|1 \hat{\;} B|=u \hat{\;} (n,m,1)=m$, and for $h=2$, $m+2$ must be divisible by 3 and thus $m-1$ cannot be a power of 2, thus the situation is a subcase of our first case above.   Note also that if $m+2$ is divisible by $h+1$, then $h \leq m-2$, since otherwise $m=h+1$, and $h+3$ cannot be divisible by $h+1$ if $h \geq 3$.  Furthermore, we have $m \geq 6$ as all smaller cases yield contradictions. 

Now observe that if $m+2$ is divisible by $h+1$, then $hm-h^2+h$ is always even; therefore, $d_0=(hm-h^2+h)/2 \in D(n)$.  But if $3 \leq h \leq m-2$ and $m \geq 6$, then $d_0 \geq m$, since 
$$2d_0-2m=(h-3)(m-h-2)+(m-6) \geq 0.$$ Therefore,     
$$u\hat{\;}(n,m,h) \leq u(n,m,h) \leq f_{d_0} = d_0 < hm-h^2 = hm-h^2+1 - (hk_2-hd+d) = |h \hat{\;} B|,$$
which is a contradiction.

{\bf Case 2}: Let us turn to the case when $B$ is a special set and $k_2 < h \leq k_1$.  We still have   
$i_{\min}=0$, but now $$i_{\max}=k_2c+(h-k_2)(c-1)=hc-h+k_2,$$ thus we see that
$$h \hat{_{i_{\min}}} B=h \hat{\;} B'$$ and $$h \hat{_{i_{\max}}} B=
\left\{ 
\begin{array}{ll}
k_2 \hat{\;} B'' + (h-k_2)\hat{\;} B'  & \mbox{if $c=1$,}\\ \\
k_2 \hat{\;} B'' + (h-k_2)\hat{\;} ((c-1)g+H) & \mbox{if $c \geq 2$.}\\
\end{array}\right.$$

We can observe right away that, if $c \geq 2$, then, by Lemma \ref{lemma any j}, $$h \hat{_{i_{\max}}} B=i_{\max} \cdot g+H;$$ in fact, in this case we have $$h \hat{_i} B=i \cdot g+H$$ for all $i_{\min}+1 \leq i \leq i_{\max}.$  Therefore, 
$$h \hat{\;} B= h \hat{\;} B' \cup \bigcup_{i=i_{\min}+1}^{i_{\max}} (i \cdot g+H),$$ and $B$ being special further implies that 
$$h \hat{\;} B'=i_{\min} \cdot g+H=i_{\max} \cdot g+H$$ and thus
$$|h \hat{\;} B|= (i_{\max}-i_{\min})\cdot d=(hc-h+k_2)d.$$

Note also that, to satisfy (\ref{g def}), $n$ must be divisible by $d_0=(hc-h+k_2)d$, and we also have
$$m=k_1+(c-1)d+k_2 < cd+k_2  \leq (hc-h+k_2)d=d_0,$$ 
and so 
$$u \hat{\;} (n,m,h) \leq u (n,m,h) \leq f_{d_0} = \left( h \left \lceil \tfrac{m}{d_0}   \right \rceil -h+1 \right)d_0=d_0=|h \hat{\;} B|.$$ 
Therefore, when $c \geq 2$, we have
$$|h \hat{\;} B| \geq u \hat{\;} (n,m,h).$$ 

Assume now that $$|h \hat{\;} B| < u \hat{\;} (n,m,h),$$ in which case $c=1$, $n$ is divisible by $(i_{\max}-i_{\min})d=k_2d$, $$g=\tfrac{n}{k_2d},$$ and we also note that $m>d$ implies that $k_2>1$.  We then have
$$
B = B' \cup B'' = \left\{ j \cdot \tfrac{n}{d} \; : \;  j=0,1,\dots,k_1-1 \right\} \cup \left\{ g+(j_0+j) \cdot \tfrac{n}{d} \; : \;  j=0,1,\dots,k_2-1 \right\},
$$
and
\begin{eqnarray*}
h \hat{\;} B & = & \bigcup_{i=0}^{k_2} (h-i) \hat{\;} B' + i \hat{\;} B'' \\ 
& = & \bigcup_{i=0}^{k_2} \left\{ i  g+ j \cdot \tfrac{n}{d} \; : \;  j=j_{\min}(i), j_{\min}(i)+1, \ldots, j_{\max}(i) \right\},
\end{eqnarray*}
where, as before,
$$j_{\min} (i)= \tfrac{(h-i)(h-i-1)}{2} + i j_0+ \tfrac{i(i-1)}{2},$$ and
$$j_{\max} (i)= (h-i)(k_1-1)-\tfrac{(h-i)(h-i-1)}{2} + i j_0+ i(k_2-1)-\tfrac{i(i-1)}{2}.$$
We find that, for $1 \leq i \leq k_2-1$, we have
\begin{eqnarray*}
j_{\max} (i)-j_{\min} (i) +1 & = & (h-i)(k_1-1)+i(k_2-1)-(h-i)(h-i-1)-i(i-1)+1 \\
& \geq & (h-i)(k_1-1)+i(k_2-1)-(h-i)(h-i-1)-i(i-1) +1+d- (k_1+k_2-1) \\
& = & (h-i-1)(k_1-h+i)+(i-1)(k_2-i)-(h-2)+d\\ 
& \geq & (h-i-1)i+(i-1)-(h-2)+d\\ 
& = & (i-1)(h-i-1)+d \\
& \geq & d, 
\end{eqnarray*}
and, therefore,
$$h \hat{\;} B=h \hat{\;} B' \cup \bigcup_{i=1}^{k_2-1} (i \cdot g+H) \cup \left(k_2 \hat{\;} B'' + (h-k_2)\hat{\;} B'\right).$$
We also find that
$$j_{\max} (0)-j_{\min} (0) +1=hk_1-h^2+1,$$ and
$$j_{\max} (k_2)-j_{\min} (k_2) +1=hk_1-h^2+1-k_2(m-2h);$$
since $B$ is special, that is, $$h \hat{_{i_{\max}}} B \subseteq h \hat{_{i_{\min}}} B,$$ we must have 
\begin{eqnarray} \label{m>=2h}
m=k_1+k_2 & \geq& 2h.
\end{eqnarray}
Therefore, we get
$$|h \hat{\;} B|=\min\{hk_1-h^2+1,d\} + (k_2-1)d.$$

Assume now that 
$$|h \hat{\;} B| < u \hat{\;} (n,m,h).$$  
Note that $d_0=k_2d \in D(n)$ and $$m=k_1+k_2 < 2d \leq k_2d=d_0,$$ and thus
$$u \hat{\;} (n,m,h) \leq u (n,m,h) \leq f_{d_0} = \left( h \left \lceil \tfrac{m}{d_0}   \right \rceil -h+1 \right)d_0=d_0.$$  Therefore, by (\ref{size of h^B}) we have
\begin{eqnarray} \label{hk_1-h^2+1 < d  }
hk_1-h^2+1 < d
\end{eqnarray}
 and thus
\begin{eqnarray} \label{|h^B| < k_2d}
|h \hat{\;} B| &=& hk_1-h^2+1 + (k_2-1) \cdot d \\
&=& hm-h^2+1-(hk_2-k_2d+d).  \label{|h^B| < hm-h^2+1 second}
\end{eqnarray}
Observe that, since
$$|h \hat{\;} B| < u \hat{\;} (n,m,h) \leq hm-h^2+1,$$ we must have  
\begin{eqnarray} \label{hk_2-k_2d+d >= 1.}
hk_2-k_2d+d \geq 1.
\end{eqnarray}

Recall that $k_1 \geq h >k_2 >1$ and $k_1+k_2>d$, so $d \leq 2k_1-2$, and thus (\ref{hk_1-h^2+1 < d  }) implies that $h \leq k_1 \leq h+2$.  We can rule out $k_1=h$ as that would contradict (\ref{m>=2h}).  If we were to have $k_1=h+2$, then from (\ref{hk_1-h^2+1 < d  }) we get $d > 2h+1$, so
$$hk_2-k_2d+d< hk_2-(k_2-1)(2h+1) = -(h+1)k_2+(2h+1) <0,$$ contradicting (\ref{hk_2-k_2d+d >= 1.}).

This leaves us with $k_1=h+1$, from which by (\ref{m>=2h}) we get $k_2=h-1$.  Since $k_2 \geq 2$, we get $h \geq 3$.  Furthermore, from (\ref{hk_1-h^2+1 < d  }) and (\ref{hk_2-k_2d+d >= 1.}) we get $h \leq 3$.    Therefore, $h=3$, $k_1=4$, $k_2=2$, $m=6$, and $d=5$.  In this case $n$ must be divisible by 10, and $g=n/10$.  We then have
$$
B' = \left\{ j \cdot \tfrac{n}{5} \; : \;  j=0,1,2,3 \right\} ,
$$
and
$$
B'' = \left\{ \tfrac{n}{10}+(j_0+j) \cdot \tfrac{n}{5} \; : \;  j=0,1 \right\},
$$
and therefore
$$
3 \hat{\;} B' = \left\{ j \cdot \tfrac{n}{5} \; : \;  j=3,4,5,6 \right\} ,
$$
and
$$
B'+ 2 \hat{\;} B'' = \left\{ (1+2j_0+j) \cdot \tfrac{n}{5} \; : \;  j=1,2,3,4\right\}.
$$
Thus, in order for $B$  to be special, we need to have $$1+2j_0 \equiv 2 \; \; \mbox{mod} \; 5,$$ thus we choose $j_0=3$.   
In summary, we have
$$B= \left\{ j \cdot \tfrac{n}{5} \; : \;  j=0,1,2,3 \right\} \cup \left\{ \tfrac{n}{10}+j \cdot \tfrac{n}{5} \; : \;  j=3,4 \right\},$$
which we recognize as the set
$$B_5 \left(n,6;4,2, \tfrac{n}{10},3 \right)$$ of the second part of Proposition \ref{special sets listed} and for which 
 $|3 \hat{\;} B|=9=hm-h^2$.
$\Box$

\section{The cases of $h=2$ and $h=3$} \label{h=2 and h=3}

Since most effort in the literature thus far has focused on the cases of $h=2$ and $h=3$, it is worth stating explicitly what we are able to using Theorem \ref{sizeofh^A_d} and Proposition \ref{special sets listed}.  

Let us start by determining $u \hat{\;} (n,m,2)$; we may assume that $3 \leq m \leq n$.  By Theorem \ref{sizeofh^A_d}, if $d \geq 3$ and $k \geq 2$, then 
$$f \hat{_d}(n,m,2)=\min\{n,f_d,2m-3\}.$$  The value of $f \hat{_d}(n,m,2)$ in the other cases can be summarized as follows.
$$\begin{array}{|c|c|c|c||c|} \hline
d & k & r & \delta_d & f \hat{_d}(n,m,2) \\ \hline \hline
1 & 1 & 1 & 0 & \min\{n,2m-3\} \\ \hline
2 & 1 & 2 & 0 & \min\{n,2m-3\} \\ \hline
2 & 2 & 2 & 1 & \min\{n,2m-4\} \\ \hline
\geq 3 & 1 & 2 & -1 & \min\{n,2m-2\} \\ \hline
\end{array}$$
Thus, we can verify that
$$u\hat{\;} (n,m,2)=\left\{
\begin{array}{ll}
\min\{u(n,m,2), 2m-4\} & \mbox{if} \; 2|n \; \mbox{and} \; 2|m; \\ \\
\min\{u(n,m,2), 2m-3\} & \mbox{otherwise.} 
\end{array}
\right.$$

Next, we use Proposition \ref{special sets listed} to determine the values of $n$ and $m$ that allow for an improvement over $u\hat{\;} (n,m,2)$.  Only one case applies: when $m-1$ is not a power of 2 and $n$ is divisible by $2m-2$, in which a special set $B$ exists with $|2\hat{\;} B|=2m-4$.  Therefore, using Theorem \ref{value of u} as well, we have the following result.  
\begin{cor} \label{cor h=2}
For all positive integers $n$ and $m$ with $3 \leq m \leq n$ we have
$$\rho\hat{\;} (\mathbb{Z}_n,m,2) \leq \left\{
\begin{array}{ll}
\min\{\rho(\mathbb{Z}_n,m,2), 2m-4\} & \mbox{if} \; 2|n \; \mbox{and} \; 2|m, \; \mbox{or} \\ 
& (2m-2)|n \; \mbox{and} \; \not \exists k \in \mathbb{N}, m=2^k+1; \\ \\
\min\{\rho(\mathbb{Z}_n,m,2), 2m-3\} & \mbox{otherwise.} 
\end{array}
\right.$$

\end{cor}
In fact, we believe that equality holds in Corollary \ref{cor h=2}; see Conjecture \ref{my conj}.

We can carry out a similar analysis for the case of $h=3$, where we may assume that $4 \leq m \leq n$.  By Theorem \ref{sizeofh^A_d}, if $d \geq 4$ and $k \geq 3$, then 
$$f \hat{_d}(n,m,3)=\min\{n,f_d,3m-8\}.$$  For the remaining cases we have the following.
$$\begin{array}{|c|c|c|c||c|} \hline
d & k & r & \delta_d & f \hat{_d}(n,m,3) \\ \hline \hline
1 & 1 & 1 & 0 & \min\{n,3m-8\} \\ \hline
2 & 1 & 1 & 0 & \min\{n,3m-8\} \\ \hline
2 & 2 & 1 & 0 & \min\{n,3m-8\} \\ \hline
3 & 1 & 3 & 0 & \min\{n,3m-8\} \\ \hline
3 & 2 & 3 & 0 & \min\{n,3m-8\} \\ \hline
3 & 3 & 3 & 2 & \min\{n,3m-10\} \\ \hline
\geq 4 & 1 & 3 & d-5 & \min\{n,3m-3-d\} \\ \hline
\geq 4 & 2 & 3 & -2 & \min\{n,3m-6\} \\ \hline
\end{array}$$
The case of $d \geq 4$ and $k=1$ is interesting: it implies that $m-1$ is divisible by $d$; thus, to minimize $3m-3-d$, we select $d=\gcd(n,m-1)$.  Note also that when $3|m$ then $\gcd(n,m-1) \neq 6$.  Therefore, 
$$u\hat{\;} (n,m,3)=\left\{
\begin{array}{ll}
\min\{u(n,m,3), 3m-3-\gcd(n,m-1)\} & \mbox{if} \; \gcd(n,m-1) \geq 8; \\ \\
\min\{u(n,m,3), 3m-10\} & \mbox{if} \; \gcd(n,m-1) = 7, \; \mbox{or} \\
& \gcd(n,m-1) \leq 5, \; 3|n, \; \mbox{and} \; 3|m; \\ \\
\min\{u(n,m,3), 3m-9\} & \mbox{if} \; \gcd(n,m-1) = 6; \\ \\
\min\{u(n,m,3), 3m-8\} & \mbox{otherwise.} 
\end{array}
\right.$$

We then examine Proposition \ref{special sets listed} for the case $h=3$ to see if we can do better.  We find two such instances: when $n$ is divisible by 10 and $m=6$, and when $n$ is divisible by $3m-9$ and $m-3$ is divisible by 5.  Observe that, in the latter case, we may assume that $m$ is even, since otherwise $d_0=(3m-9)/2 \in D(n)$ and thus
$$u\hat{\;} (n,m,3) \leq u(n,m,3) \leq f_{d_0} =d_0 \leq 3m-10.$$  Therefore, we have the following.

\begin{cor} \label{cor h=3}
For all positive integers $n$ and $m$ with $4 \leq m \leq n$ we have
$$\rho\hat{\;} (\mathbb{Z}_n,m,3) \leq \left\{
\begin{array}{ll}
\min\{u(n,m,3), 3m-3-\gcd(n,m-1)\} & \mbox{if} \; \gcd(n,m-1) \geq 8; \\ \\
\min\{u(n,m,3), 3m-10\} & \mbox{if} \; \gcd(n,m-1) = 7, \; \mbox{or} \\
&  \gcd(n,m-1) \leq 5, \; 3|n, \; \mbox{and} \; 3|m, \; \mbox{or} \\
&  \gcd(n,m-1) \leq 5, \; (3m-9)|n, \; \mbox{and} \; 5|(m-3); \\ \\
\min\{u(n,m,3), 3m-9\} & \mbox{if} \; \gcd(n,m-1) = 6, \; \mbox{or} \\
& m=6 \; \mbox{and} \; 10|n \; \mbox{but} \; 3 \not | n; \\ \\
\min\{u(n,m,3), 3m-8\} & \mbox{otherwise.} 
\end{array}
\right.$$

\end{cor}

We have performed a computer search for all $m$-subsets of $\mathbb{Z}_n$ with $4 \leq m \leq n \leq 40$, and in each case we found that equality holds in Corollary \ref{cor h=3}.

{\bf Acknowledgments.}  I wish to thank Anna Llad\'o and Alain Plagne for their friendly exchanges regarding my counterexamples for their conjectures, and Gyula K\'arolyi for helpful comments on the manuscript.  I am also grateful for Marcin Malec for generating valuable computational data.


\begin{thebibliography}{99}

\bibitem{AloNatRuz:1995a} N. Alon, M. B. Nathanson, and I. Ruzsa, Adding Distinct Congruence Classes Modulo a Prime, {\em Amer. Math. Monthly}, {\bf 102} (1995) 250--255.

\bibitem{AloNatRuz:1996a} N. Alon, M. B. Nathanson, and I. Ruzsa, The Polynomial Method and Restricted Sums of Congruence Classes, {\em J. Number Theory}, {\bf 56} (1996) 404--417.


\bibitem{Cau:1813a} A.-L. Cauchy, Recherches sur les nombres, {\em J. \'Ecole Polytechnique} {\bf 9} (1813) 99--123.

\bibitem{Dav:1935a} H. Davenport, On the addition of residue classes, {\em J. London Math. Soc.} {\bf 10} (1935) 30--32.

\bibitem{Dav:1947a} H. Davenport, A historical note, {\em J. London Math. Soc.} {\bf 22} (1947) 100--101.


\bibitem{DiaHam:1994a} J. A. Dias da Silva and Y. O. Hamidoune, Cyclic space for Grassmann derivatives and additive theory.  {\em Bull. London Math. Soc.}, {\bf 26} (1994) 140--146.

\bibitem{EliKer:1998a} S. Eliahou and M. Kervaire, Sumsets in Vector Spaces over Finite Fields, {\em J. Number Theory}, {\bf 71} (1998) 12--39.

\bibitem{EliKer:2001a} S. Eliahou and M. Kervaire, Restricted Sumsets in Finite Vector Spaces: The Case $p=3$, {\em Integers}, {\bf 1} (2001) \#A02.

\bibitem{EliKer:2001b} S. Eliahou and M. Kervaire, Restricted sums of sets of cardinality $1+p$ in a vector space over ${\bf F}_p$, {\em Discrete Math.}, {\bf 235} (2001) 199--213.

\bibitem{EliKer:2007a} S. Eliahou and M. Kervaire, Some extensions of the Cauchy--Davenport Theorem, {\em Electronic Notes in Discrete Math.}, {\bf 28} (2007) 557--564.

\bibitem{EliKerPla:2003a} S. Eliahou, M. Kervaire, and A. Plagne, Optimally small sumsets in finite abelian groups, {\em J. Number Theory}, {\bf 101} (2003) 338--348.

\bibitem{ErdGra:1980a} P. Erd\H{o}s and R. L. Graham, {\em Old and New Problems and Results in Combinatorial Number Theory}. L'Enseignement Math\'ematique, Geneva (1980).

\bibitem{ErdHei:1964a} P. Erd\H{o}s and H. Heilbronn, On the addition of residue classes (mod $p$).  {\em Acta Arith.}, {\bf 9} (1964) 149--159.

\bibitem{GalGre:2002a} L. Gallardo, G. Grekos, et al., Restricted addition in $\mathbb{Z}/n\mathbb{Z}$ and an application to the Erd\H{o}s--Ginzburg--Ziv problem.  {\em J. London Math. Soc. (2)}, {\bf 65} (2002) 513--523.

\bibitem{GirGriHam:2012a} B. Girard, S. Griffiths, and Y. O. Hamidoune, $k$-sums in abelian groups.  {\em Combin. Probab. Comput.}, {\bf 21} (2012) no. 4, 582--596.

\bibitem{HamLlaSer:2000a} Y. O. Hamidoune, A. S. Llad\'o, and O. Serra, On Restricted Sums.  {\em Combin. Probab. Comput.}, {\bf 9} (2000) 513--518.


\bibitem{Kar:2003a} Gy. K\'arolyi, On restricted set addition in abelian groups.  {\em Ann. Univ. Sci. Budapest, E\"otv\"os Sect. Math.}, {\bf 46} (2003) 47--54.

\bibitem{Kar:2004a} Gy. K\'arolyi, The Erd\H{o}s--Heilbronn problem in abelian groups.  {\em Israel J. Math.}, {\bf 139} (2004) 349--359.

\bibitem{Kar:2005a} Gy. K\'arolyi, An inverse theorem for the restricted set addition in abelian groups.  {\em J. Algebra}, {\bf 290} (2004) 557--593.

\bibitem{Kar:2006a} Gy. K\'arolyi, A note on the Hopf--Stiefel function.  {\em European J. Combin.}, {\bf 27} (2006) 1135--1137.


\bibitem{Kne:1956a} M. Kneser, Summenmengen in lokalkompakten abelschen Gruppen, {\em Math. Zeitschr.}, {\bf 66} (1956) 88--110.


\bibitem{Lev:2000a} V. F. Lev, Restricted set addition in groups I: The classical setting, {\em J. London Math. Soc. (2)} {\bf 62} (2000) 27--40.

\bibitem{Lev:2002a} V. F. Lev, Three-fold Restricted Set Addition in Groups, {\em European J. Combin.} {\bf 23} (2002) 613--617.

\bibitem{Lev:2005a} V. F. Lev, Restricted set addition in Abelian groups: results and conjectures, {\em J. Théor. Nombres Bordeaux} {\bf 17} (2005) 181--193.

\bibitem{Nat:1996a} M.~B. Nathanson, {\em Additive Number Theory: Inverse Problems and the Geometry of Sumsets}.  Springer-Verlag New York, 1996.

\bibitem{Pla:2003a} A. Plagne, Additive number theory sheds extra light on the Hopf--Stiefel $\circ$ function, {\em Enseign. Math., II S\'er}, {\bf 49}:1--2 (2003) 109--116.

\bibitem{Pla:2006a} A. Plagne, Optimally small sumsets in groups, I. The supersmall sumset property, the $\mu_G^{(k)}$ and the $\nu_G^{(k)}$ functions, {\em Unif. Distrib. Theory}, {\bf 1}:1 (2006) 27--44.

\bibitem{Pla:2006b} A. Plagne, Optimally small sumsets in groups, II. The hypersmall sumset property and restricted addition, {\em Unif. Distrib. Theory}, {\bf 1}:1 (2006) 111--124.

\bibitem{Sha:1984a} D. Shapiro, Products of sums of squares, {\em Expo. Math.}, {\bf 2} (1984) 235--261.

\bibitem{Wan:2008a} G. Wang, On restricted sumsets in abelian groups of odd order, {\em Integers}, {\bf 8} (2008) \#A22.



\end{thebibliography}
\end{document}